\documentclass[12pt]{article}

\usepackage{verbatim}
\usepackage{mathpazo}
\usepackage{graphicx}
\usepackage{amsfonts}
\usepackage{amssymb}
\usepackage{amsmath}
\usepackage{url}

\usepackage{lmodern}
\usepackage[a4paper]{geometry}
\usepackage{amsthm}
\usepackage{amsmath,amssymb,amsfonts}

\theoremstyle{plain}
\newtheorem{thm}{Theorem}
\newtheorem{prop}[thm]{Proposition}

\newtheorem{lem}[thm]{Lemma}

\newtheorem{fact}[thm]{Fact}

\newcommand{\R}{\ensuremath{\mathbb{R}}} 

\newcommand{\be}[1]{\begin{equation}\label{#1}}
\newcommand{\ee}{\end{equation}}
\newcommand{\D}{\displaystyle} 

\theoremstyle{remark}

\newtheorem{rmq}{Remark}

\def\var{\varepsilon}

\title{More on functional and quantitative versions of the isoperimetric inequality}
\author{Erik Thomas}
\begin{document}

\maketitle

\begin{abstract}
The goal of the present paper is to discuss new functional extensions of the isoperimetric inequality, and quantitative versions involving the Wasserstein's distances. 
\end{abstract}

\section{Introduction}
We shall work on the Euclidean space $(\R^n,\cdot, |\cdot|)$. The sharp (anisotropic) isoperimetric inequality  can be stated as follow: given a convex body $K\subset \R^n$ (having zero in  its interior), if we denote by 
$$n'= \frac{n}{n-1}$$
the Lebesgue conjugate to $n$,  we have for every Borel set $E\subset\R^n$ that
\be{classicalisop}
p_K \left( E \right)   \ge  n \left| K \right|^{\frac{1}{n}} \left| E \right|^{\frac1{n'} }  ,
\ee
with equality if $ E= \lambda K + a$ for some $\lambda > 0$ and $a \in \R^n.$ 
Here 

$$p_K(E)= \liminf_{\var \to 0} \frac{|E + \var K| - |E|}{\var}.$$ 
Equivalently, if $E$ has a regular enough boundary $\partial E$, then 
$$p_K(E) = \int_{\partial E} h_K(-\nu(x)) \, d\mathcal H^{n-1} (x) ,$$
where 
$$h_K(z) := \sup_{y\in K} y \cdot z, \qquad \forall z \in \R^n.$$
is the support function of the body $K$, $\nu(x)$ is the outer unit normalto $\partial E$ at $x\in \partial E$, and $\mathcal H^{n-1}$ stands for the $(n-1)$-dimensional Hausdorff measure on $\partial E$.  The classical Euclidean isoperimetric inequality corresponds to the case when $K=B_2^n=\{|\cdot|\le 1\}$, the Euclidean unit ball. 

We want to analyse functional versions of~\eqref{classicalisop}. Replacing $E$ by a locally Lipschitz function $f:\R^n \to \R_+$ is  standard. We have decided to work with nonnegative functions recalling that for $f$ with values in $\R$ we can apply the result to $\left| f \right|$ and use the fact that for $f \in W^{1,1}_{\mathrm{loc}} \left( \R^n \right),$ we have, almost-everywhere, $\nabla \left| f  \right|= \pm \nabla f$. The functionnal inequality takes the same form
\be{classicalisop2}
p_K \left( f \right)  \ge n \left| K \right|^{\frac{1}{n}}\left(\int_{\R^n} f^{n'}  \right)^{1/n'},
\ee
with equality if $f= 1_E$ (provided the gradient term below is understood as a capacity of the bounded variation function $1_E$). Here, 
$$p_K(f) = \int_{\R^n} h_K(-\nabla f (x)) \, dx .$$ 
The inequality~\eqref{classicalisop2} can be proven directly using a mass transportation method, as observed by Gromov, see the appendix of \cite{ms}. In the case of the Euclidean ball, $K=B_2^n$, we recover
$$p_{B_2^n}(f) =\big\| |\nabla f |\big\|_{L^1(\R^n)}.$$

Extending the convex body $K$ to a (convex) function or measure is less obvious. First, one needs to have a proper extension of the notion of support function $h_K$ for a convex function $V$. Actually, the integral term $\int h_K(-\nabla f)$ needs a proper interpretation, so that non only a convex function will enter the game, but also some "convex measure" (in the terminology of Borell \cite{bor1} and \cite{bor2}) associated to it.  
This has been studied recently in several papers. In particular, in~\cite{kl} corresponding extensions of the isoperimetric inequality~\eqref{classicalisop2} have been proposed. See also~\cite{cofr} and~\cite{mr}. Here we will establish a new inequality that has the advantage to contain the geometric versions~\eqref{classicalisop} and~\eqref{classicalisop2}. We will do that by picking a good category of convex measures. First, we need to introduce some notation. Let $V:\R^n \to \R\cup\{+\infty\}$ be a nonnegative convex function such that $\D Z_V:= \int_{\R^n} \left( 1+\frac1{n-1} V(x) \right)^{-n}  d x < +\infty$. We associate to $V$ the probability measure 
 \be{defmuV}
 d\mu_V(x) = \frac{1}{Z_V} \,   \big(1+\frac1{n-1} V(x) \big)^{-n}\, dx .
 \ee
 Our  generalization for $p_K \left( f \right)$ is as follows, for $f : \R^n \to \R_+$ locally Lipschitz we put
 
$$p_{V} \left( f \right) := \int_{\R^n}  V^{\ast} \left( \frac{- \nabla f}{f^{\frac{n}{n-1}}} \right) f^{\frac{n}{n-1}} +  \left(\int_{\R^n} V  d \mu_V \right)  \left(\int_{\R^n} f^{n'}  \right)$$
where $V^\ast$  is the Legendre's transform of $V$, 
$$V^{\ast}  \left( y \right) = \sup_{x \in \R^n} x \cdot y - V \left( x \right), \qquad \forall y \in \R^n.$$
In particular, we have following inequality (known as Young's inequality): 
\begin{equation}
\forall x, y \in \R^n, \qquad x \cdot y \le V \left( x \right) + V^{\ast} \left( y \right),
\label{young1}
\end{equation}
with equality when $y= \nabla V \left( x \right)$. Note that when 
$$V = 1_K^\infty:=\begin{cases}
0 & \textrm{ on } K\\
+\infty & \textrm{ outside } K\
\end{cases}$$
is the "indicatrix" of a convex set $K$,  then $p_V \left(f \right) =   \int_{\R^n} h_K \left( - \nabla f \right)  = p_K \left(f \right)$ since $V=0$ $\mu_V$ almost-everywhere and $V^\ast = h_K$. The general isoperimetric-Sobolev inequality, we can get is then as follow.

\begin{thm}\label{theo:iso}
Let $V$ be a nonnegative convex function with $Z_V=\int_{\R^n} (1+ \frac{V}{n-1})^{-n} < +\infty$ and $\mu_V$ the associated probability measure~\eqref{defmuV}. Then, for every nonnegative locally Lipschitz function $f$ on $\R^n$ we have

\begin{equation}
p_{V} \left( f \right)  \ge \left[ n \, Z_V^{\frac1n} {\small \int_{\R^n}  (  1 + \frac{1}{n-1} V )  d \mu_V} \right] \, \| f \|_{L^{n'} \left( \R^n \right)},
\label{thm1}
\end{equation}
and, when $V$ is finite, with equality when $f \left( x \right) =  \left( 1 + \frac{1}{n-1} V \left( x - a \right) \right)^{- \left( n - 1 \right)}$ with $a \in \R^n.$
\end{thm}
Thanks to the remark prior to the theorem, we see that when $V=1_K^\infty$, inequality~\eqref{thm1} becomes exactly~\eqref{classicalisop2}.

The second topic of the present paper is mainly independent of what we discussed so far, although  based again on mass transport methods. We aim at presenting some quantitative forms of the geometric isoperimetric inequality~\eqref{classicalisop} that involve a Kantorovich-Rubinstein (or Wassertein) distance cost to an extremizer.  For $u : \mathbb{R}^n \to \mathbb{R}^n$ a Borel map and  $\mu$ a measure in $\mathbb{R}^n$, we write $u_{\sharp} \mu$ for the measure defined by
$$ u_{\sharp} \mu \left( M \right) := \mu \left( u^{-1} \left( M \right) \right) $$ for all Borel sets $M \subseteq \mathbb{R}^n$. It is called the push-forward of $\mu$ through $u$.

 For $\mu$ and $\nu$ two probability measures in $\R^n$, and a (cost) function $c:\R^n \times \R^n \to \R_+$, 
 we define the Kantorovich-Rubinstein or Wasserstein transportation cost $ \mathcal W_c \left( \mu , \nu \right)$ by
\begin{eqnarray*}
\mathcal W_c \left( \mu , \nu \right)&=&\inf_{T : \mathbb{R}^n \to \mathbb{R}^n : T_{\sharp} \mu=\nu    }  \int_{\mathbb{R}^n}  c(x,T \left( x \right))  d\mu \left( x \right) \\
 &=& \inf_\pi \iint_{\R^n \times \R^n} c(x,y)  d\pi(x,y)
\end{eqnarray*}
 where the infimum is taken over  probability measures $\pi$ on $\R^n\times \R^n$ that have $\mu$ and $\nu$ as marginals, respectively. 
 For $1 \le p \le +\infty$ and $c(x,y) = |y-x|^p$, the $p$-th power classical $p$-Kantorovich-Rubinstein distance,  $W_p^p$, is recovered. We refer to~\cite{vi} for details.

Our measures will be uniform measures on the sets $E$ and $K$. Given a Borel set  $E$, we will denote by $\lambda_E$ the Lebesgue measure restricted to $E$ and normalized to be a probability measure,
$$ d \lambda_E(x) =\frac{1_E \left( x \right)}{ \left| E \right|} dx$$
  where $1_E$ is the indicator function of the set $E$.  Given a Borel set $E$, we will denote by $\widetilde{E}$ the homothetic of volume one of the set $E$, namely  
$$\widetilde{E} = \frac{1}{ \left| E \right|^{\frac{1}{n}} } E .$$
Note that for  $u \in GL_n \left( \mathbb{R} \right)$ we have $\lambda_{u \left( E\right)}=  u_{\sharp} \lambda_E$. And with some abuse of notation,  we denote for $t>0$ by $t_\sharp \mu$ the image of $\mu$ under the dilation by $t$, we have
$$ t _\sharp \lambda_E = \lambda_{tE} \quad \textrm{and}\quad \lambda_{\widetilde{E}}=\frac{1}{|  E |^{\frac1n}}\, \mbox{}_{\sharp} \, \lambda_E.$$

Our cost function will depend on the set $E$. Recall that for a probability measure $\mu$ on $\R^n$, its \emph{Cheeger constant} $ D_{\mathrm{Che}} \left( \mu \right)$ is  the best (i.e. largest)  constant such that the following inequality holds for all Borel sets $A$:
$$\mu^{+} \left( A \right) \ge D_{\mathrm{Che}} \left( \mu \right)  \min \left\{ \mu \left( A \right) , 1 - \mu \left( A \right) \right\}$$ where $\mu^+$ denotes the measure of the perimeter (or Minkowski content) associated to $\mu.$  It can be defined by:
$$\mu^+ \left( A \right) := \liminf_{\epsilon \to 0} \frac{  \mu \left( A_{\epsilon}  \right)  -\mu \left( A \right) }{\epsilon},$$ where $A_{\epsilon}=\left\{ x \in \R^n : \mathrm{dist} \left( x , A \right) < \epsilon   \right\}.$ Equivalently, if we denote by $h_{p,q}(\mu)$ the best nonnegative constant for which the inequality 
$$ \left( \int \left| \nabla f \left( x \right) \right|^q  d\mu(x) \right)^{\frac1q} \ge h_{p,q}(\mu) \left(  \int \left| f\left( x \right)  - \int f\, d \mu  \right|^p   d\mu(x) \right)^{\frac1p}$$ 
holds for all $f \in W^{1,1}_{\mathrm{loc}} \cap L^p \left( \mu \right)$, then $h_{1,1}(\mu) \le D_{\mathrm{Che}} (\mu) \le 2 h_{1,1}(\mu)$.
Let $\mathcal{F}$ be the convex, increasing function defined on $\mathbb{R}_+$ by 
$$\mathcal{F} \left( t \right):=t -\log(1+t).$$ 
The function $\mathcal F$  behaves like $t^2$ for $t$ small and like $t$ for $t$ large, and
$$\min  \left\{ t^2, t \right\} \le \mathcal F(t) \le 2 \min \left\{ t^2, t  \right\}, \qquad \forall t \ge 0.$$
This function appears in several mass transport proofs to give a remainder term, to instance in \cite{fmp}, \cite{bk}, \cite{ceg}.

Given a probability measure $\mu$, we will use the following cost $c:\R^n \times \R^n \to \R_+$ that is also used in \cite{ce}:
\be{defcost}
c_\mu(x,y) := \mathcal F \left( D_{\mathrm{Che}} \left( \mu \right) \, |y-x| \right)
\ee
which behaves like $D_{\mathrm{Che}}(\mu)^2 \, |y-x|^2$ for small distances, and like $D_{\mathrm{Che}}(\mu) \, |y-x| $ for large ones.

Our main result is the following extension of the isoperimetric inequality.

\begin{thm}\label{quantitative}
Let $K$ be a convex body on $\R^n$.  Given a Borel set $E\subset \R^n$ with locally Lipschitz boundary and $\int_{\widetilde{E}} x dx = \int_{\widetilde{K}} x dx$, we have
\begin{equation}
R(E,K):= \frac{ p_K \left( E \right)   }{ n \left| K \right|^{\frac{1}{n}} \left| E \right|^{\frac{n}{n-1}}  }   - 1  
\ge 
\frac{c}{n}  \mathcal W_{c_{\lambda_{ \widetilde{E}}}} \big(\lambda_{  \widetilde{E} } ,   \lambda_{  \widetilde{K} } \big),
\label{bound1}
\end{equation}
for some universal constant $c>0$, and as a consequence
\begin{equation}
R(E,K) \ge \frac{c}{n} \mathcal{F} \left( D_{\mathrm{Che}}(\lambda_E ) W_1 (\lambda_{  {E} } ,   \lambda_{  {K} } )  \right)
\label{bound1bis}
\end{equation}
\end{thm}

We emphasize here a weakness of this result: the remainder term depends on $E$ (on the Cheeger constant of $\widetilde{ E}$, precisely). But in some geometric problems, $\widetilde{E}$ will not be too wild :  it will belong to a family of sets for which we have a good control on $D_{\mathrm{Che}}( \lambda_{\widetilde{E}})$, as we will see later. 
Since the condition  $\int_{\widetilde{E}} x dx = \int_{\widetilde{K}} x dx$  can always be achieved by translating $E$, we can drop this assumption provided the transportation term is replaced by 
$\min_{v\in \R^n} \mathcal{W}_{c_{\lambda_{\widetilde E}}} (\tau_v \lambda_{\widetilde E}, \lambda_{\widetilde K})$ where $\tau_v \nu$ is the image of the measure $\nu$ by the translation by $v$ in $\R^n$.  Note that we still have equality if $E= \lambda K$ for some $\lambda > 0.$ 

When $E$ is convex, it is known that $D_{\mathrm{Che}}(\lambda_E)>0 $. In this case, we also know that up to numerical constants, $D_{\mathrm{Che}}(\lambda_E)$ is the same as $h_{2,2} \left( \lambda_E \right)$, the Poincar\'{e} constant associated to $E$ (or the inverse of the spectral gap).

Let us compare our results to existing quantitative Sobolev and isoperimetric inequalities, obtained by Figalli-Maggi-Pratelli. In \cite{fmp}, there is a quantitative isoperimetric inequality (the numerical constant we use are the improved ones obtained by Segal \cite{se}):
\begin{equation}
p_K \left( E \right)    \ge n \left| K \right|^{\frac{1}{n}}  |E|^{ \frac{1}{n'} }  \left( 1   + \frac{C}{ n^7}  A_K \left( E \right)^2 \right)   ,
\label{fmp}
\end{equation}
where $A_K \left( E \right):=\inf \left\{  \frac{ \left| E \Delta \left(  x_0 + r K \right) \right|  }{\left| E \right|} : x_0 \in \R^n, r^n \left| K \right| = \left| E \right|  \right\}$ and $C$ is a numerical constant.

This result of  Figalli-Maggi-Pratelli  is much deeper and in general stronger than ours, since it is universal (the bound does not depend on geometry of $E$, as in our case). We can note however that  the quantity $ \frac{C}{n^7} A_K \left( E \right)^2 $ decreases to $0$ when the dimension $n$ goes to $+\infty.$

Actually, there are some particular cases in which our result might give a better bound, both in fixed dimension and when the dimension grows. The reason is that the transportation cost term can be rather large. For instance, we will give examples where our remainder, $\mathcal{F} \left( D_{\mathrm{Che}}(\lambda_E )  W_1 (\lambda_{  {E} } ,   \lambda_{  {K} } )  \right)$ decreases slower than $\frac{1}{n^7}$ of the inequality~\eqref{fmp}.

The rest of the paper is organized as follows. In the next section, we collect some results on optimal transportation theory. Then, we will prove our two Theorems above.
In a final section, we will compute our reminder term in several situation of interest arising in convex geometry. 
\newline
\newline
I would like to thank my Professor Dario Cordero-Erausquin for his encouragements, his careful reviews and his many useful discussions.

\section{Proof of Theorem 1}

We first give some background about optimal transportation.

\subsection{Background on optimal transportation}

The following Theorem, due to Brenier \cite{be} and refined then by McCann \cite{mc1}, is the main result in optimal transportation.

\begin{thm}

If $\mu$ and $\nu$ are two probability measures on $\mathbb{R}^n$ and $\mu$ absolutely continuous with respect to Lebesgue measure, then there exists a convex function $\phi$ such that $T=\nabla \phi$ transports $\mu$ onto $\nu.$ Moreover, $T$ is uniquely determined $\mu$ almost-everywhere. 
\end{thm}
That means that for every nonnegative Borel function $b : \mathbb{R}^n \to \mathbb{R}_+,$
\begin{equation}
\int_{ \mathbb{R}^n} b \left( y \right)  d \nu \left( y \right) = \int_{ \mathbb{R}^n} b \left( T \left( x \right) \right) d \mu \left( x \right)  .
\label{cosette}
\end{equation}
If $\mu$ and $\nu$ have densities, say $F$ and $G$,~\eqref{cosette} becomes

\begin{equation}
\int_{ \mathbb{R}^n} b \left( y \right) G \left( y \right)  d y= \int_{ \mathbb{R}^n} b \left( \nabla \phi \left( x \right) \right) F \left( x \right)  d x  .
\label{fantine}
\end{equation}
If $\phi$ is $C^2$ the change of variables $y=\nabla \phi (x)$ in~\eqref{fantine} gives the Monge-Ampère equation, for $F \left( x \right) d x $ almost-every $x \in \R^n$:

\begin{equation}
F \left( x \right)= G \left(  \nabla \phi \left( x \right) \right) \mathrm{det} \left( D^2 \phi \left( x \right) \right),
\label{valjean}
\end{equation}
where $D^2 \phi$ is the hessian matrix of $\phi.$

\begin{rmq}

When $T$ is the Brenier map between $\lambda_E$ and $\lambda_K$ with $E$ and $K$ two convex bodies with same volume,~\eqref{valjean} is simpler:

$$\mathrm{det} \left( D^2 \phi \left( x \right) \right) = 1, $$
for $\lambda_E$ almost-every $x \in E.$

\end{rmq}

The question of regularity of $\phi$ can be asked because, in the previous equality~\eqref{valjean}, $\phi$ seemed to be required $C^2.$ In fact, this is not the case, as it was established by McCann \cite{mc2} that we can give an almost-everywhere sense to~\eqref{valjean} by rather standard arguments from measure theory. This almost-everywhere theory is sufficient for most applications, including the one in the present paper but it requires some further arguments that will be discussed later. 

\subsection{Proof of Theorem 1: the inequality}

Let us recall the frame. Let $V:\R^n \to \R_+ \cup \left\{ +\infty \right\}$ a nonnegative convex function such that $Z_V=\int_{\R^n} \left( 1 + \frac{1}{n-1} V \left( x \right) \right)^{-n}  d x < + \infty.$ So we define the probability measure $\mu_V$ by $$d \mu_V \left( x \right) = \frac{1}{Z_V} \left( 1 + \frac{1}{n-1} V \left( x \right) \right)^{-n} d x.$$ 
Let $f: \R^n \to \R_+$ a Borel function such that $0 < \int_{\R^n} f^{\frac{n}{n-1}} < + \infty.$ So we can define the probability measure $\mu$ by

$$d \mu \left( x \right)= \frac{f^{\frac{n}{n-1}}   \left( x \right) }{  \int_{\R^n} f^{\frac{n}{n-1}} } d x.$$
Let $T=\nabla \varphi$ the Brenier between $\mu$ and $\mu_V$ and we start by studying the regularity of $\varphi.$ It is sufficient to prove the Theorem for measures $\mu$ and $\mu_V$ whose support is $\R^n.$ We also can assume that $f$ is the convolution of a function compactly supported and a mollifier, so that $f$ is smooth and converges rapidly to $0$  at $+\infty.$ Then, it is known that prove that $\varphi \in W^{2,1}_{\mathrm{loc}} \left( \R^n \right)$ and the following equality

\begin{equation}
\int_{\R^n} f \, \Delta \varphi = - \int_{\R^n} \nabla f \cdot \nabla \varphi
\label{ibp1}
\end{equation}
is valid. Let us prove now the first part of Theorem 1.

\begin{proof}

We first need the following Fact.

\begin{fact}{\cite{cenv}} Let $d \mu \left( x \right)=F \left( x \right) d x$ and $d \nu \left( y \right)=G \left( y \right) d y$ two probability measures on $\R^n.$ Let $T=\nabla \varphi$ the Brenier map between $\mu$ and $\nu.$ Then, the following inequality holds:

$$\int_{\R^n}  G^{1 - \frac1n} \le \frac1n \int_{\R^n} F^{1 - \frac1n} \Delta \varphi .$$ 

\end{fact}

Let us give the proof this Fact for completeness. 

\begin{proof}

We start with Monge-Amp\`{e}re equation, for $\mu$ almost-every $x\in \R^n,$ we have:

$$F\left( x \right)=G \left( \nabla \varphi \left( x \right) \right) \det \left( D^2 \varphi \left( x \right) \right).$$
Then, for $\mu$ almost-every $x \in \R^n$ and thanks to arithmetic-geometric inequality:

\begin{equation}
G^{-\frac1n} \left(  \nabla \varphi \left( x \right) \right) \le F^{-  \frac1n} \left( x \right) \frac{\Delta \varphi \left( x \right)}{n}.
\label{faitfonda} 
\end{equation}
An integration with respect to $d \mu \left( x \right) = F \left( x \right) d x$ gives:

\begin{eqnarray*}
\frac 1n \int_{\R^n} F^{1- \frac1n} \left( x \right) \Delta \varphi \left( x \right) d x &\ge& \int_{\R^n} G^{- \frac1n} \left(  \nabla \varphi \left( x \right)  \right) d x \\
&\underbrace{=}_{~\eqref{fantine}}& \int_{\R^n} G^{1 - \frac1n}  \left( x \right) d x.
\end{eqnarray*}

\end{proof}
If we apply this Fact to our situation, it gives:

$$ \int_{\R^n}  \frac{ \left(  1 + \frac{1}{n- 1} V \right)^{-n + 1} }{ Z_V^{\frac{1}{n'}} } \le \frac1n \int_{\R^n }  \frac{f}{ \| f \|_{L^{n'} \left(  \R^n \right)} } \Delta \varphi,$$
and

$$n \, \| f \|_{L^{n'} \left( \R^n \right)} \, \int_{\R^n} \frac{ \left(  1 + \frac{1}{n-1} V \right)^{-n+1}  }{Z_V^{\frac{1}{n'}}} \le \int_{\R^n} f \Delta \varphi .$$
As $n \, \| f \|_{L^{n'} \left( \R^n \right)} \, \int_{\R^n} \frac{ \left(  1 + \frac{1}{n-1} V \right)^{-n+1}  }{Z_V^{\frac{1}{n'}}} = n \, \| f \|_{L^{n'} \left( \R^n \right)} Z_V^{\frac1n} \, \int_{\R^n}  \left( 1 + \frac{1}{n-1} V \right) d \mu_V,$ we now have the following lines:

\begin{eqnarray*}
n  \| f \|_{L^{n'} \left(  \R^n \right)} Z_V^{\frac1n} \int_{\R^n}  \left( 1 + \frac{1}{n-1} V \right) d \mu_V  &\le& \int_{\R^n} f  \, \Delta \varphi \\
&\underbrace{=}_{\eqref{ibp1}}& \int_{\R^n} \left( - \nabla f \right) \cdot \nabla \varphi \\
&=& \int_{\R^n} \left(  \frac{- \nabla f}{f^{\frac{n}{n-1}}}  \cdot  \nabla \varphi \right) f^{\frac{n}{n-1}} \\
&\underbrace{\le}_{\eqref{young1}}& \int_{\R^n} V^\ast \left(  \frac{- \nabla f}{f^{\frac{n}{n-1}}} \right) f^{\frac{n}{n-1}} + \int_{\R^n} V \left( \nabla \varphi \right) f^{\frac{n}{n-1}} \\
&=& \int_{\R^n} V^\ast \left(  \frac{- \nabla f}{f^{\frac{n}{n-1}}} \right) f^{\frac{n}{n-1}} + \left( \int_{\R^n}  V \circ T  d \mu \right) \left( \int_{\R^n} f^{n'} \right) \\
&=& \int_{\R^n} V^\ast \left(  \frac{- \nabla f}{f^{\frac{n}{n-1}}} \right) f^{\frac{n}{n-1}} + \left( \int_{\R^n} V d \mu_V \right)  \left( \int_{\R^n} f^{n'} \right).
\end{eqnarray*}

\end{proof}

\subsection{Case of equality}

In this subsection, we establish that the inequality~\eqref{thm1} becomes an equality when $f \left( x \right) =  \left( 1 + \frac{1}{n-1} V \left( x  \right) \right)^{- \left( n - 1 \right)}$ with $V: \R^n \to \R$ a finite convex function.  Note that in this case the Brenier map $T=\nabla \varphi$ is $T \left( x \right) = x$ so $D^2 \varphi = I$ and Monge-Amp\`{e}re equation is, for $\mu$ almost-every $x \in \R^n,$

\begin{equation}
f^{\frac{n}{n-1}} \left( x \right)= \left( 1 + \frac{1}{n-1}  V \left( T \left( x \right) \right) \right)^{-n} \underbrace{ \det \left(  D^2 \varphi \left( x \right) \right) }_{=1}.
\label{mabis}
\end{equation}
If we come back to the proof of the inequality~\eqref{thm1}, we remark that we use only two inequalities: the inequality in Fact 4 (which is an arithmetic-geometric inequality) and Young's inequality.

We note that in our case case, the inequality in Fact 4 in an equality since $D^2 \varphi=I$ so $\det \left(  D^2 \varphi \left( x \right) \right)=\frac{\Delta \varphi \left( x \right)}{n}.$ \\
Let us treat now Young's inequality. We have an equality in Young's inequality if (and only if) for $\mu$ almost-every $x \in \R^n$

\begin{equation}
\frac{- \nabla f \left( x \right)}{f^{\frac{n}{n-1}} \left( x \right) }=\nabla V \left( T \left( x \right) \right)= \nabla V(x).
\label{young}
\end{equation}
To get this equality, let us take the $-\frac1n$ power in~\eqref{mabis} to get

$$f^{\frac{1}{n-1}} \left( x \right)=1  +  \frac{1}{n-1} V \left( T \left( x \right) \right)=1 + \frac{1}{n-1} V \left( x \right).$$
If we compute the gradient of the previous line, we find~\eqref{young}.

\begin{rmq}

One can prove that~\eqref{thm1} is an equality \emph{if and only if} $f \left( x \right) =  \left( 1 + \frac{1}{n-1} V \left( x - a \right) \right)^{- \left( n - 1 \right)}$ with $a \in \R^n.$ We decided not to prove this because it is technical. Let us speak about this. If we have an equality in~\eqref{thm1}, we have an equality in the inequality in Fact 4 and an equality in Young's inequality.

An equality in Fact 4 means that $\det \left(  D^2 \varphi \left( x \right) \right) = \frac{\Delta \varphi \left( x \right)}{n}$ so the matrix $D^2 \varphi \left( x \right)$ has only one eigenvalue, say $\lambda \left( x \right)$ and $D^2 \varphi \left( x \right) = \lambda \left( x \right) I.$ The main difficulty is to show that the function $\lambda $ is constant, which is the case when $\varphi$ is $C^2$ smooth (the details are analyzed in~\cite{cent}). If we assume that, it is easy to conclude that, up to translations, $\mu=\mu_V.$

\end{rmq}

\section{Proof of Theorem 2}

Here we establish Theorem 2. It was noted by Figalli, Maggi and Pratelli~\cite{fmp} (and Segal~\cite{se}) that for this kind of result, the general situation follows from the case the two bodies have same volume, equal to one. For completeness, let us recall the argument.
Let $E$ Borel set and $K$ a convex body in $\mathbb{R}^n.$ The following Lemma establishes a link between $ R \left( E , K \right)=\frac{p_K \left(  E \right)}{ n \left| K
\right|^{\frac{1}{n}} \left| E \right|^{\frac{1}{n'}} } - 1 $ and $R \left( \widetilde{E} , \widetilde{K}  \right)=\frac{p_ {\widetilde{K}} \left(  \widetilde{E} \right)}{ n \left| \widetilde{K}
\right|^{\frac{1}{n}} \left| \widetilde{E} \right|^{\frac{1}{n'}} } - 1$ where $\widetilde{E}$ and $\widetilde{K}$ are respectively  $\frac{E}{\left| E  \right|^{\frac{1}{n}}}$ et $\frac{K}{\left| K  \right|^{\frac{1}{n}}}.$

\begin{lem}
With the previous notations, 
$$R \left( E , K \right)=\frac{p_K \left(  E \right)}{ n \left| K
\right|^{\frac{1}{n}} \left| E \right|^{\frac{1}{n'}} } - 1 =  \frac{p_ {\widetilde{K}} \left(  \widetilde{E} \right)}{ n \left| \widetilde{K}
\right|^{\frac{1}{n}} \left| \widetilde{E} \right|^{\frac{1}{n'}} } - 1=R \left( \widetilde{E} , \widetilde{K}  \right).$$ 
\end{lem}

\begin{proof}
Let us note that for all $\epsilon > 0$, we have
$$ \frac{ \left| E +\epsilon K \right| - \left| E \right| }{\epsilon } = \left| E \right|^{\frac{1}{n}} \left| K \right|^{\frac{1}{n'}} \frac{\left| \widetilde{E} +\epsilon \frac{ \left| K \right|^{\frac{1}{n}} }{\left| E \right|^{\frac{1}{n}}} \widetilde{K} \right| - \left| \widetilde{E} \right|  }{\epsilon\frac{ \left| K \right|^{\frac{1}{n}} }{\left| E \right|^{\frac{1}{n}}}}.$$
By taking the limit, we get the equality.
\end{proof}
Therefore, if we have established Theorem 2 for two sets of volume one, we have the general statement by applying it to $\widetilde{E}$ and $\widetilde{K}$. So in the rest of this section, $E$ is a Borel set with smooth boundary and $K$ a convex body, both with volume one, $|E|=|K|=1$.  

As in~\cite{fmp} and~\cite{se}, the argument to establish Theorem 2 starts with optimal transportation. The following Lemma gives a first minimization for the deficit $R \left( E , K \right)=\frac{ p_K \left(  E \right) }{ n \left| K \right|^{\frac{1}{n}} \left| E \right|^{\frac{1}{n'}} }  - 1$.

\begin{lem}{\cite{fmp}} Let $E$ and $K$ two convex bodies in $\mathbb{R}^n$ with same measure $1.$ Let $T=\nabla \phi$ the Brenier map between the measures $\lambda_E$ and $\lambda_K.$ We note by $0 < \lambda_1 \le \cdots \le\lambda_n$ the eigenvalues of $D^2 \phi.$
Then, we have the following inequality
\begin{equation}
R \left( E , K \right) \ge \int_{\mathbb{R}^n} \left( \lambda_A - \lambda_G \right) d \lambda_E, 
\label{hernani}
\end{equation}
where $\lambda_A=\frac{\lambda_1 + \cdots + \lambda_n}{n}$ and $\lambda_G=\Pi_{i=1}^n \lambda_i^{\frac{1}{n}}.$

\end{lem}

Before briefly recalling the proof of this Lemma, let us speak about the regularity of the optimal transport. It is known, see~\cite{ca}, that when $T=\nabla \phi$ is the brenier map between $d \mu \left( x \right)=f \left( x \right)  d x$ and $d \nu \left( y \right) = g \left( y \right) d y$ two probability measures supported on two open bounded sets, respectively $E$ and $K,$ with $f$ and $g$ are $\alpha$-H\"{o}lder, bounded and with $\frac{1}{f}$ and $\frac{1}{g}$ bounded too, then $\phi \in C^{2,\beta} \left( E \right)$ for all $0 < \beta < \alpha.$ 

\begin{proof}

As $\left| E \right|=\left| K \right|=1,$ we can write
$$n \left| K \right|^{\frac{1}{n}} \left| E \right|^{\frac{1}{n'}}= \int_{\mathbb{R}^n}
n \left(  \mathrm{det} \left(   D^2 \phi    \right) \right)^{\frac{1}{n}}  d \lambda_E =\int_E n \left(  \mathrm{det} \left( D^2 \phi \right) \right)^{\frac{1}{n}},$$ because $\mathrm{det} \left( D^2 \phi \right)=1$ thanks to the Remark 1. The arithmetic-geometric inequality gives
$$ n \left| K \right|^{\frac{1}{n}} \left| E \right|^{\frac{1}{n'}} \le
\int_{E} \mathrm{div} \, T  \left( x \right)  d x.$$
The divergence theorem provides
$$ \int_{E} \mathrm{div} \, T  \left( x \right)   d x = \int_{\partial
E}  T \left( x \right) \cdot \nu_E \left( x \right) 
 d \mathcal{H}^{n-1} \left( x \right) ,$$
where $\mathcal{H}^{n-1}$ is the $\left(n-1\right)$-dimensional Haussdorf measure. By definition of the support function $h_K,$ of $K$, and since $T(E)\subseteq K$ we therefore get  
$$   n \left| K \right|^{\frac{1}{n}} \left| E \right|^{\frac{1}{n'}} \le
\int_{E} \mathrm{div} T  \left( x \right)   d x \le \int_{\partial
E}  h_K \left(  \nu_E \left( x \right) \right) d \mathcal{H}^{n-1}
\left( x \right)=p_K \left(  E \right) .$$
Thus

$$\frac{ p_K \left(  E \right) }{n \left| K \right|^{\frac{1}{n'}} \left| E \right|^{\frac{1}{n}} } -  1 
\ge \int_{\mathbb{R}^n} \left( \frac{ \mathrm{div} T }{n} - 1 \right)  d \lambda_E 
\ge \int_{\mathbb{R}^n} \left( \lambda_A - \lambda_G \right)  d \lambda_E .$$

\end{proof}
To go on, we need a quantitative version of the arithmetic-geometric inequality. The following result is due to Alzer~\cite{a}.

\begin{lem}{\cite{a}} Let $0< \lambda_1 \le \cdots \le \lambda_n.$ Let $\lambda_A=\frac{\lambda_1 + \cdots + \lambda_n}{n}$ and $\lambda_G=\Pi_{i=1}^n \lambda_i^{\frac{1}{n}}.$ We have
\begin{equation}
\sum_{i=1}^n \left(  \lambda_i - \lambda_G \right)^2 \le 2 n \lambda_n \left(  \lambda_A - \lambda_G \right).
\label{hugo}
\end{equation}

\end{lem}

We can now complete the proof of Theorem 2. If $T=\nabla\varphi$ is the Brenier map between $\lambda_E$ and $\lambda_K$ for two bodies $E$ and $K$ of volume $1$. With the previous notations, if we use~\eqref{hernani} and~\eqref{hugo}, we get
$$ \left(  \lambda_A - \lambda_G \right) \ge \frac{1}{2n} \frac{ \| D^2 \phi - \mathrm{Id} \|_{\mathrm{HS}}^2 }{\lambda_n} \ge \frac{1}{2n} \frac{\| D^2 \phi - \mathrm{Id} \|_{\mathrm{HS}}^2}{1+ \| D^2 \phi - \mathrm{Id} \|_{\mathrm{HS}}}  \ge  \frac{c}{n}  \mathrm{tr} \left( \mathcal F \left(  \|  D^2 \theta  \|_{\mathrm{HS}} \right) \right) , $$
where $\theta \left( x \right)= \phi \left( x \right) - \frac{\left| x \right|^2}{2}$ and $\|\cdot\|_{\mathrm{HS}}$ refers to the Hilbert-Schmidt norm of a $n\times n$ matrix.  Let us remark that $\lambda_G=1,$ thanks to, once again, Remark 1. So we have 

\begin{equation}
R \left( E , K \right) \ge \frac{c}{n} \int_{\mathbb{R}^n} \mathrm{tr} \left( \mathcal F \left(  \|  D^2 \theta  \|_{\mathrm{HS}} \right) \right).
\label{remainder}
\end{equation}
The treatment of this term is stated in the next Lemma and we refer to~\cite{ce}.

\begin{lem}{\cite{ce}} Let $\mu$ a probability measure on $\R^n$ absolutely continuous with respect to the Lebesgue measure and $\theta \in W_{\mathrm{loc}}^{2,1} \left( \R^n \right)$ with $D^2 \theta + Id \ge 0$ almost-everywhere. We assume $\left| \nabla \theta \right| \in L^1 \left( \mu \right)$ and $\int_{\R^n} \nabla \theta  d \mu =0.$ Then,

$$ \int_{\R^n}  \mathrm{tr} \left( \mathcal F \left( D^2 \theta \right) \right) d \mu \ge c \int_{\R^n}  \mathcal F \left( D_{\mathrm{Che}} \left( \mu \right) \left| \nabla \theta \right| \right) d \mu,$$
for some numerical constant $c>0.$

\end{lem}
Note that our assumption $\int_{E} x dx = \int_{K} x dx$ rewrites as $\int_{E} \nabla \theta = 0,$ so if we use the previous Lemma with $\mu=\lambda_E,$ in~\eqref{remainder} we find

\begin{eqnarray*}
R \left( E , K \right) &\ge& \frac{c}{n} \int_{\R^n} \mathcal F \left( D_{\text{Che}} \left|  \nabla \theta  \right| \right) d \lambda_E \\
&\ge& \frac{c}{n} \mathcal W_{c_{\lambda_E}} \left(  \lambda_E, \lambda_K \right).
\end{eqnarray*}

\section{Some examples}

Here we give some examples where our result (i.e. Theorem 2) gives good bounds for the remainder term, better thant the one in~\cite{fmp}. We will give an example in dimension 2 where our remainder term, depending on a parameter, can be as large as we want and an example in dimension $n$. We recall that the remainder term in~\eqref{fmp} is bounded by 1 when $E$ and $K$ have for measure 1 and decreases to $0$ with $\frac{1}{n^7}$ when the dimension $n$ grows.

\subsection{In dimension 2}

In this section, we give a toy example in dimension 2. Let, for $\alpha > 0,$ $E_{\alpha} = \left[ -\frac{\alpha}{2} , \frac{\alpha}{2} \right] \times \left[ - \frac{1}{2 \alpha} , \frac{1}{2 \alpha} \right] $ and $K_{\alpha} = \left[ - \frac{\alpha^2}{2} , \frac{\alpha^2}{2} \right] \times \left[ - \frac{1}{2 \alpha^2} , \frac{1}{2 \alpha^2} \right].$ We will prove the following:

\begin{prop}

With the previous notations, we have:
$$\lim_{\alpha \to +\infty} D_{\mathrm{Che}} \left( \lambda_{E_{\alpha}} \right) W_1 \left( \lambda_{E_{\alpha}} , \lambda_{K_{\alpha}} \right)= +\infty.$$

\end{prop}

\begin{proof}
As $W_1 \left( \lambda_{E_{\alpha}} , \lambda_{K_{\alpha}} \right) = W_1 \left(   \lambda_{K_{\alpha}} , \lambda_{E_{\alpha}} \right),$ we give an estimation of the last term. Let $T=\nabla \phi$ the Brenier map which transports the measure $\lambda_{K_{\alpha}}$ onto the measure $\lambda_{E_{\alpha}}$ (by Monge-Amp\`{e}re equation, it verifies $\det D^2 \phi = 1,$ in particular it preserves the volume). Let $K_{\alpha}'= \left[ \frac{\alpha^2}{4} , \frac{\alpha^2}{2} \right] \times \left[ - \frac{1}{2 \alpha^2} , \frac{1}{2 \alpha^2} \right].$ Then, we have:

$$ W_1 \left( \lambda_{E_{\alpha}} , \lambda_{K_{\alpha}} \right) = \int_{\R^n} \left| T \left( x \right) - x \right|   d \lambda_{K_{\alpha}} \left( x \right)  \ge \int_{K_{\alpha}'} \left| T \left( x \right) - x \right| d x \ge \left| K_{\alpha}' \right| \mathrm{dist} \left( K_{\alpha}' , E_{\alpha} \right).$$
So, we have

$$W_1 \left( \lambda_{E_{\alpha}} , \lambda_{K_{\alpha}} \right) \ge \frac{1}{4} \left( \frac{\alpha^2}{4} - \frac{\alpha}{2} \right).$$

\end{proof}

We need an estimation of $D_{\mathrm{Che}} \left( \lambda_{E_{\alpha}} \right).$ This constant could be computed explicitly but it is rather easier to compare this constant to the Poincar\'{e} constant  $h_{2,2} \left( \lambda_{E_{\alpha}} \right).$ Indeed, it is known that, up to numerical constants, that $D_{\mathrm{Che}} \left( \lambda_{E_{\alpha}} \right) $ is the same as $h_{2,2} \left( \lambda_{E_{\alpha}} \right)$, see~\cite{le} and~\cite{em}. As $ \lambda_{E_{\alpha}}=\lambda_{ \left[ -\frac{\alpha}{2} , \frac{\alpha}{2} \right] } \otimes \lambda_{ \left[ - \frac{1}{2 \alpha} , \frac{1}{2 \alpha} \right]},$ then $h_{2,2} \left( E_{\alpha} \right) \ge \min \left\{ h_{2,2} \left( \lambda_{ \left[ -\frac{\alpha}{2} , \frac{\alpha}{2} \right] }\right) , h_{2,2} \left(   \lambda_{ \left[ - \frac{1}{2 \alpha} , \frac{1}{2 \alpha} \right] }  \right)  \right\},$ see~\cite{bh} and~\cite{bo1}. Since $h_{2,2} \left( \lambda_{  \left[ - a , a \right]  } \right)=\frac{\pi}{a}$ for $a>0,$,it follows with Theorem 2 that:

$$R \left( E_{\alpha} , K_{\alpha} \right) \ge \frac{1}{4} \mathcal{F} \left(   \frac{c \pi}{2 \alpha} \left( \frac{\alpha^2}{4} - \frac{\lambda}{2}  \right)  \right) ,$$
for some numerical constant $c>0$. In particular the remainder term $R \left( E_{\alpha} , K_{\alpha} \right)$ is not bounded when $\alpha$ grows to $+\infty$ whereas the remainder in~\eqref{fmp} remains bounded.

\subsection{Estimation of $W_1 \left( \lambda_K , \lambda_L\right)$ for $K$ and $L$ isotropic convex bodies}

Here $K$ and $L$ are two convex bodies with measure $1.$ We say that a convex body $K$ is in isotropic position if $\left| K \right|=1,$ it is centered and there exists $\alpha > 0,$ such that

$$\int_K \left| x \cdot y \right|^2 \,  dx = \alpha \left| y \right|^2, \qquad \forall y \in \R^n.$$
For an isotropic convex body $K$, we define its isotropic constant $L_K (=\sqrt \alpha)$ by:

$$L_K^2 =\frac 1n \int_K \left| x \right|^2 dx.$$
We also define for any convex isotropic body $K$

$$M \left( K \right) = \frac{1}{\sqrt{n}} \int_K \left| x \right| d x.$$
Using H\"{o}lder inequality and Borell deviation inequality \cite{bor1, bor2}, we have:

$$c L \left( K \right) \le M \left( K \right) \le L \left( K \right),$$
for some numerical constant $c>0.$ For backgrounds, we refer to~\cite{bgvv}. Our goal is here is to prove the following Proposition.

\begin{prop}

Let $K$ and $L$ two convex bodies of volume $1$ in isotropic position. Then, the following estimation for $W_1 \left( \lambda_K , \lambda_L \right)$ holds:

\begin{equation}
\sqrt{n} \left| M \left( K \right) - M \left( L \right) \right| \le W_1 \left( \lambda_K , \lambda_L \right) \le c \left( L_K + L_L \right) \sqrt{n} + 8,
\label{inequa}
\end{equation}
for some numerical constant $c>0.$ In connection with some isoperimetric estimates, we are mainly interested with lower bounds.

\end{prop}

We are mainly interested in the lower bound provided by this Proposition. The upper bound (which is far from sharp when $K$ are $L$ are closed to each other) is stated only to emphasize that the generic expected order of magnitude is $\sqrt n$ on both sides.

\begin{proof}
We first work on the left hand-side of~\eqref{inequa}. Let us recall the dual of $W_1 \left( \lambda_K , \lambda_L \right),$ known as Kantorovich's duality:

\begin{equation}
W_1 \left( \mu , \nu \right) = \sup_{\phi \, 1-\mathrm{Lip}} \left\{ \int_{\R^n} \phi \, d  \mu - \int_{\R^n} \phi \, d \nu  \right\}.
\label{kanto_utile}
\end{equation}
If we take in~\eqref{kanto_utile}, $ \phi \left( x \right) = \left| x \right| $ or $- \left| x \right|,$ we have:

\begin{eqnarray*}
W_1 \left( \lambda_K , \lambda_L \right) &\ge& \left| \int_{\R^n} \left| x \right| d \lambda_K \left( x \right) -  \int_{\R^n} \left| x \right| d \lambda_L \left( x \right)   \right| \\
&=&\left| M \left( K \right) - M \left( L \right)  \right|.
\end{eqnarray*}

Let us treat now the right hand-side of~\eqref{inequa}. Let $T$ be a transport map (in particular, it verifies $\det \left( \nabla T \right) = 1$) which transports the measure $\lambda_K$ onto the measure $\lambda_L,$ so $W_1 \left( \lambda_K , \lambda_L \right) \le \int_{\R^n} \left| T \left( x \right) - x \right|  d \lambda_K \left( x \right) = \int_K \left| T \left( x \right) - x \right|  d x.$ Let us  recall a deep result of Paouris, see~\cite{p}.

\begin{thm}{\cite{p}} There exists a numerical constant $c> 0$ such that if $K$ is an isotropic convex body in $\R^n,$ then

\begin{equation}
 \left|  \left\{  x \in K : \left| x \right| \ge c \sqrt{n} L_K t \right\} \right| \le \exp \left( - \sqrt{n} t \right), \qquad \forall t \ge 1.
\label{paou}
\end{equation}

\end{thm}

Let $t \ge 1$ such that $\exp \left( - \sqrt{n} t  \right) \le \frac1n$ (note that $t=1$ works, we will set this value for $t$) so the following sets $K_1= \left\{  x \in K : \left| x \right| < c \sqrt{n} L_K \right\}$ and $L_1= \left\{  x \in L : \left| x \right| < c \sqrt{n} L_L \right\}$ have their volumes bigger than $1- \frac1n.$ Finally, let $K_2=T^{-1} \left( L_1 \right).$ Since $\det \left( \nabla T \right)=1,$ we have $\left| K_2 \right| = \left| L_1 \right|.$ We can now conclude thanks to the following inequality:

\begin{equation}
W_1 \left( \lambda_K , \lambda_L \right) \le \int_{K} \left| T \left( x \right) - x \right|  d x = \int_{K_1 \cap K_2} \left| T \left( x \right) - x \right|  d x + \int_{K \backslash \left( K_1 \cap K_2 \right)} \left| T \left( x \right) - x \right|  d x.
\label{equaa}
\end{equation} 
Let us estimate the two remaining integrals. Thanks to the definitions of $K_1$ and $K_2,$ we have: $ \int_{K_1 \cap K_2} \left| T \left( x \right) - x \right|  d x \le c \left( L_K + L_L \right) \sqrt{n}.$ For the second, we need the fact that $K,L \subseteq B \left( 0 , \sqrt{n \left( n + 2 \right)} \right) \subseteq B \left( 0 ,  2n \right),$ see~\cite{kls}. Then, for $x \in K \backslash \left( K_1 \cap K_2 \right),$ we have $\left| T \left( x \right) - x \right| \le 4n$ and $\left| K \backslash \left( K_1 \cap K_2 \right) \right| \le \frac{2}{n},$ that gives $\int_{K \backslash \left( K_1 \cap K_2 \right)} \left| T \left( x \right) - x \right|  d x \le 8.$ Going back to~\eqref{equaa}, we finally have:

$$W_1 \left( \lambda_K , \lambda_L \right) \le \int_{K} \left| T \left( x \right) - x \right|  d x  \le c \left( L_K + L_L \right) \sqrt{n} + 8.$$

\end{proof}

\vskip1cm 
\noindent Erik Thomas \\
  Institut de Math\'ematiques de Jussieu, \\ Universit\'e Pierre et Marie Curie - Paris 6, \\
  75252 Paris Cedex 05, France \\
 \verb?erik.thomas@imj-prg.fr?

\end{document}